\documentclass[11pt,twoside]{amsart}
\usepackage{hyperref, tikz, fullpage}
\usepackage{placeins, enumerate,longtable}
\usepackage{amsthm, amssymb,amsmath,latexsym}
\usepackage{mathtools}
\usepackage[utf8]{inputenc}
\usepackage[T1]{fontenc}
\usepackage{enumitem,color}
        
\usepackage{caption}
\usepackage{lscape}
\newtheorem{theorem}{Theorem}[section]

\theoremstyle{definition}
\newtheorem{definition}[theorem]{Definition}
\newtheorem{example}[theorem]{Example}

\numberwithin{equation}{section}
\begin{document}
\title{Two generation of finite simple groups}
\author{Yash Arora}
\email{arorayash2015@gmail.com}
\author{Anupam Singh}
\email{anupamk18@gmail.com}
\address{Indian Institute of Science Education and Research (IISER) Pune,  Dr. Homi Bhabha Road Pashan, Pune 411008, India}
\thanks{This article is based on the MS thesis of the first named author submitted to IISER Pune. The second named author acknowledges support of the DST-RFBR project via grant number INT/RUS/RFBR/P-288.}

\begin{abstract}
This expository article revolves around the question to find short presentations of finite simple groups. This subject is one of the most active research areas of group theory in recent times. We bring together several known results on two-generation and $(2,3)$-generation of finite simple groups and how it impacts computational group theory.
\end{abstract}
\maketitle
\section{Introduction}
The groups usually arise from symmetries of an object. One of the ways groups arise naturally is while studying topological invariants, e.g., the fundamental group, homotopy groups etc. The groups, in this situation, naturally come with generators and relations. Often the groups are realised as a quotient of certain infinite groups, namely the free groups. This gives rise to a completely new way of looking at groups (as opposed to the definition and examples given to us in our undergraduate classes), and is studied under the subject of combinatorial group theory and geometric group theory. This is quite different from the usual notion of groups, where we know all of the elements, and how to multiply them. 

In the modern times, this way of defining groups has gained more importance due to the subject of computational group theory (see~\cite{heo}). Due to advancements in computational power, it is natural to expect that one can use computers to solve various mathematical problems. There are excellent mathematical packages such as GAP, MAGMA, SAGEmath, to name a few, where one can work with groups. Usually several well known groups, for example, symmetric groups, finite simple groups, matrix groups, etc. are implemented in these packages. Each of these packages allows one to do further computations within those groups. Thus from the point of view of implementation, it is not efficient to define a group on the computer using all of its elements along with its multiplication table. Thus, a ``small'' presentation is a better way to implement a group. This brings us to two different points of view of looking at these packages: how to develop better algorithms to implement these groups, and another, how to use these packages to work with these groups and compute further within them. In this article we keep our focus on the family of finite simple groups, more specifically, the finite classical groups.
\begin{definition}
A non-trivial group $G$ is said to be simple if it has no non-trivial proper normal subgroup, i.e., the only normal subgroups of $G$ are $\{e\}$ and $G$.
\end{definition} 
\noindent While studying finite groups it is natural to ask if we can classify all finite simple groups. The answer is in affirmative and we briefly recall the classification of finite simple groups (CFSG), and refer the reader to an excellent book by Wilson~\cite{wi} on this subject. The classification of finite simple group is one of the main achievements of the last century. We all should aspire to know at least the statement. The finite simple groups FSG (see~\cite{wi} page 3 for details) are broadly in four families: 
\begin{enumerate}
\item[FSG1:] Cyclic groups of prime order $\mathbb Z/p\mathbb Z$ (these are the only Abelian groups). 
\item[FSG2:] Alternating groups $A_n$ for $n\geq 5$. 
\item[FSG3:] Finite groups of Lie type: classical types and exceptional types. 
\begin{enumerate}
\item[FSG3a:] Classical groups of Lie type
\begin{enumerate}
\item[$A_l$ type:] The projective special linear groups, $PSL_{l+1}(q)$, for $l\geq 1$ except $PSL_2(2)$ and $PSL_2(3)$.
\item[$B_l$ type:] The projective quotient of the commutator of orthogonal groups, $P\Omega_{2l+1}(q)$, for $l\geq 3$ and $q$ odd.
\item[$C_l$ type:] The projective symplectic groups, $PSp_{2l}(q)$, for $l\geq 2$ except $PSp_4(2)$.
\item[$D_l$ type:] The projective quotient of commutator of orthogonal groups, $P\Omega_{2l}^{+}(q)$, for $l\geq 4$.
\item[${}^2A_l$ type:] The projective special unitary groups, $PSU_{l+1}(q)$, for $l\geq 2$ except $PSU_3(2)$.
\item[${}^2D_l$ type:] The projective quotient of commutator of orthogonal groups, $P\Omega_{2l}^{-}(q)$, for $l\geq 4$. 
\end{enumerate}
\item[FSG3b:] Exceptional groups of Lie type:
\begin{enumerate}
\item $G_2(q), q\geq 3$; $F_4(q)$; $E_6(q)$; $E_7(q)$; $E_8(q)$; ${}^2E_6(q)$; ${}^3D_4(q)$.
\item ${}^2B_2(2^{2n+1}), n\geq 1$; ${}^2G_2(3^{2n+1}), n\geq 1$; ${}^2F_4(2^{2n+1}), n\geq 1$; ${}^2F_4(2)'$.
\end{enumerate}
\end{enumerate}
\item[FSG4:] The $26$ sporadic groups named as follows:
\begin{enumerate}
{\setlength\itemindent{100pt} \item[Mathieu groups:] $M_{11}, M_{12}, M_{22}, M_{23}, M_{24}$;
\item[Leech lattice groups:] $Co_1, Co_2, Co_3, McL, HS, Suz, J_2$;
\item[Fischer groups:] $Fi_{22}, Fi_{23}, Fi_{24}'$;
\item[Monster groups:] $\mathbb M, \mathbb B, Th, HN, He$;
\item[Pariahs:] $J_1, J_2, J_3, O'N, Ly, Ru$. }
\end{enumerate}
\end{enumerate}
In all of the above $q$ is a prime power $p^a$, which indicates the size of underlying finite field $\mathbb F_q$. There is a very small number of duplication in the above list. 

In this article, by simple groups we always mean non-abelian simple groups. We begin with the basic idea of groups given by generators and relations. A broad general problem is to determine which of the (finite) groups are generated by $2$ of its elements. Such groups would be a quotient of the free group on two generators $\mathcal F_2$. However, this could be a wild problem, thus one restricts to understand which of the finite simple groups are two-generated. This has been very well studied over the years, and reasonably good answers are known. While working with finite simple groups, it was noted that many of these groups can be generated by two elements, one of order $2$ and another of order $3$. It's a fundamental theorem in the subject, due to Feit and Thompson, that every FSG has an element of order $2$.  Thus, the question is to further determine, which of the finite simple groups are $(2,3)$-generated, i.e., generated by an element of order $2$ and another of order $3$. This problem can be also thought of as determining the quotients of the free product group $C_2\star C_3$.  In this article we briefly present some of the work done in this direction.

\subsection*{Acknowledgments}  
The authors would like to thank Professor B. Sury for his lecture on this topic in the ``Workshop on Group Theory 2019'' held at IISER Pune. We also thank Dr Uday Bhaskar Sharma and Professor Marco Antonio Pellegrini for their feedback on this article. We thank the referee(s) for suggestions which improved the readability of this paper.

\section{Generators and relations for a group}
Let $G$ be a group. 
\begin{definition}
A presentation of the group $G$ is
\begin{equation}\label{presentation} 
G=\langle S \mid R\rangle,
\end{equation}
where $S$ is a subset of $G$ which generates $G$, and $R$ is a set of words on $S$ called relations, i.e., $G\cong F(S)/N(R)$ where $F(S)$ is the free group on $S$ and $N(R)$ is the normal subgroup of $F(S)$ generated by the set of relations $R$. 
\end{definition}
\noindent The presentation~\ref{presentation} is said to be a finite presentation if both $S$ and $R$ are finite. In this paper, we discuss only finite presentations, even though the groups may be finite or infinite. Let us begin with recalling some examples of presentations. 

\begin{example}
The symmetric group $S_n$ has a Coxeter presentation given by,
 $$S_n= \langle s_1, s_2, \ldots, s_{n-1} \mid s_i^2, (s_is_{i+1})^3, (s_is_j)^2,\  1\leq i<j\leq n-1, |i-j|\geq 2\rangle.$$ 
 Here we can identify $s_i$ with the transposition $(i,i+1)$ to get the isomorphism. 
\end{example}
\noindent The reflection groups are defined to be certain subgroups of the orthogonal group $O_n(\mathbb R)$ generated by some order $2$ elements. There is a more general theory of Coxeter groups, and we refer an interested reader to the book by Humphreys~\cite{hu} on this subject.  

\begin{example}
$S_n=\langle (1,2), (1,2,\cdots, n)\rangle$. Thus, $S_n$ is a quotient of $C_2\star C_n$. In fact, the presentation with respect to this generating set is as follows:
$$S_n = \left \langle x, y \mid x^2, y^n, (xy)^{n-1}, (xy^{-1}xy)^3, (xy^{-j}xy^j)^2, 2\leq j \leq \lfloor n/2\rfloor \right\rangle.$$
\end{example}

\begin{example}[Dihedral group]
The group of symmetries of a regular $n$-gon in the plane, is the Dihedral group $D_n$ with $2n$ elements. Its presentations are as follows: 
$$D_n=\langle r,s \mid r^n, s^2, (sr)^2\rangle = \langle s_1, s_2 \mid s_1^2, s_2^2, (s_1s_2)^n \rangle.$$
The infinite dihedral group is 
$$D_{\infty}= C_2\star C_2=\langle x,y \mid x^2, y^2\rangle.$$ 
It has another presentation $D_{\infty}=\langle r, s \mid s^2, (rs)^2 \rangle$. Using this, one can easily see that the finite dihedral group $D_n$ is a quotient of $D_{\infty}$.
\end{example}

\begin{example}[Modular group]\label{modular-group}
The group $PSL_2(\mathbb Z)=SL_2(\mathbb Z)/\{\pm I \}$ is called the modular group, where
$$SL_2(\mathbb Z)=\left\{ \begin{pmatrix} a&b \\ c&d \end{pmatrix} \mid a,b,c,d \in \mathbb Z, ad-bc=1 \right\}.$$
It is an infinite group, and has the presentation: $PSL_2(\mathbb Z)\cong C_2\star C_3=\langle x,y \mid x^2, y^3\rangle$. We refer the reader to an article by Conrad~\cite{co} for the proof of this interesting fact. 
\end{example}

\begin{example}[Hurwitz groups]\label{hurwitz}
The group $\Delta=\langle x,y \mid x^2, y^3, (xy)^7\rangle$ is called a $(2,3,7)$-triangular group, and any finite quotient of this group is called a Hurwitz group. These groups arise in the study of automorphisms of Riemann surfaces. We refer the reader to an article by Conder~\cite{co1} for further study. 
\end{example}

\subsection{Word problem}
When a group $G$ is given by generators $S$ and relations $R$, we can write a random element of the group $G$ as a word in the generators. However, several words might represent the same group element. Thus, it is an important problem to decide when these words represent the same group element, which also amounts to finding words that represent the identity. If we have an algorithm to decide if a word represents the identity of the group $G$, then we say that we have a solution to the word problem in $G$. 

\subsection{Cayley graph}
Let $G$ be a group and $S$ a generating set of $G$. We assume $1\notin S$. The Cayley graph $\Gamma(G,S)$, of the group $G$ with respect to $S$, is defined as follows. 
\begin{definition}
The vertex set $V(\Gamma)$ for the graph $\Gamma(G,S)$ is given by the elements of the group $G$. The edge set $E(\Gamma) = \{(g, sg) \mid g\in G, s\in S\}$, i.e., there is an edge from $g_1$ to $g_2$ if $g_2=sg_1$ for some $s\in S$. 
\end{definition}
\noindent This graph is usually directed. However, when the set $S$ is symmetric (i.e. $S=S^{-1}$) the graph is undirected. The Cayley graph of finite simple groups provide examples of expander family of graphs (see~\cite{kln,bggt,bg}) thus playing an important role in this subject. In this article we will not go into this aspect, instead we refer an interested reader to the book by Tao~\cite{ta}.

\subsection{Computational group theory}
When a group is defined by generators and relations on the computer, it becomes  important to have a fast algorithm which produces a random element of the group. Often groups are defined as a subgroup of the symmetric group or matrix groups, as these groups are easier to implement. We refer the reader to look at the book~\cite{heo} on this subject. This subject has given rise to several interesting research problems and associated projects to solve these problems such as ``group recognition project''.

\section{Two generation problem for finite simple groups}
Let $G$ be a group. Do there exist two elements in $G$, such that the group $G$ is generated by those two elements? A further question is if we can put restriction on the order of elements, for example, can we have one of these elements of order $2$. 
\begin{definition}
A group $G$ is said to be two generated if it has two elements $r,s$ in $G$ such that $G=\langle r, s \rangle$. 
\end{definition}
\noindent A two generated group is a quotient of the free group on two generators $\mathcal F_2$. Classifying two generated finite groups would be a very general problem (for all groups), thus a restricted question is, to determine which of the finite simple groups are two generated. This problem is often referred to as the ``two generation problem''.  We mostly focus on the family of classical groups (FSG3a in our notation).

\subsection{Chevalley-Steinberg generators for the groups of Lie type}
One of the largest family of finite simple groups, is of the groups of Lie type. Chevalley (see~\cite{ch}), and his work extended by Steinberg (see~\cite{st2}), gave a uniform method, starting from simple Lie algebras over $\mathbb C$, to construct these groups over any field, by providing an explicit set of generators. We briefly explain this process here, and refer to the book by Carter~\cite{ca} on this subject for further details. 

Let $\mathfrak L$ be a simple Lie algebra over $\mathbb C$. Let $\Phi$ be the corresponding reduced root system of $\mathfrak L$. Chevalley proved that there exists a basis of $\mathfrak L$ such that all the structure constants are integers. Such a basis is called a Chevalley basis, and it essentially means that $\mathfrak L$ can be defined over $\mathbb Z$. Let $k$ be a field. For each $r\in \Phi$ and $t\in k$, there are certain automorphisms $x_r(t)$ of the Lie algebra $\mathfrak L$. Let $G$ be the group generated by these elements $x_r(t)$. These generators are called Chevalley generators of the corresponding group. Using the adjoint representation, Chevalley and Steinberg proved that the groups thus obtained are simple groups, and over the finite field $k=\mathbb F_q$, they give the family of finite simple groups (FSG3 in our notation). Let us understand this process through some examples.
\begin{example}\label{sl-gen}
In the group $SL_n(\mathbb F_q)$, the elements $x_{i,j}(t)=I+ te_{i,j}$ for $1\leq i \neq j \leq n$, where $t\in \mathbb F_q$ and $e_{i,j}$ is the matrix with $1$ at the $ij^{th}$ place and $0$ elsewhere, are called Chevalley generators. The Gaussian elimination algorithm, using row-column operations, provides an algorithmic proof that these elements generate the group $SL_n(\mathbb F_q)$. In fact, a smaller subset of this set,
$$\{x_{i,i+1}(t), x_{i+1, i}(s) \mid t,s \in \mathbb F_q, 1\leq i \leq n-1\},$$
generates the group $SL_n(\mathbb F_q)$. This is because we can get all of the remaining Chevalley generators by taking commutators of these ones.
\end{example}
\begin{example}\label{sp-gen}
Let $q$ be odd. Let $J=\begin{pmatrix} 0 & I_l \\-I_l & 0 \end{pmatrix}$, and $Sp_{2l}(\mathbb F_q)=\{X\in GL_{2l}(\mathbb F_q) \mid  {}^tXJX=J\}$ be the symplectic group. Following~\cite{ca}, we use the index set for the matrix of size $2l$ as $1, \ldots, l, -1, \ldots, -l$.  The group $Sp_{2l}(\mathbb F_q)$ is generated by the following set of Chevalley generators: 
$$\{ x_{i,j}(t), i\neq j; x_{i,-j}(t), i< j; x_{-i,j}(t), i< j, x_{i,-i}(t), x_{-i,i}(t) \mid t\in \mathbb F_q\},$$ 
where $x_{i,j}(t)=I+t(e_{i,j}-e_{-j,-i}); x_{i,-j}(t)=I+t(e_{i,-j}+e_{j,-i}); x_{-i,j}(t)=I+t(e_{-i,j}+ e_{-j,i}); x_{i,-i}(t)= I+te_{i,-i};\text { and } x_{-i,i}(t)=I+te_{-i,i}$.
 One can further compute that the simple generators (corresponding to simple roots) with their negative counterparts are enough to generate $Sp_{2l}(\mathbb F_q)$. The remaining generators can be produced using the commutators of these. For example, if we work with $Sp_6(\mathbb F_q)$, the set 
 $$\{x_{1,2}(t), x_{2,3}(t), x_{3, -1}(t), x_{2,1}(t), x_{3,2}(t), x_{-1, 3}(t) \mid t\in \mathbb F_q\}$$ 
 generates this group. 
\end{example}

We refer the reader to~\cite{ms,bmss} for analogue of Gaussian elimination algorithm in orthogonal groups, symplectic groups and unitary groups, which provides an algorithmic  proof of generation of the corresponding groups using their Chevalley generators. Notice that the number of Chevalley generators is usually large, as it varies with the field size and the matrix size. However, if we restrict to the base field $\mathbb F_p$, this is still an interesting generating set.

\begin{example}
The group $SL_2(\mathbb F_p)$ is generated by the two elements, $\begin{pmatrix} 1& 1 \\ 0 &1\end{pmatrix}$ and $\begin{pmatrix} 1& 0 \\ 1&1 \end{pmatrix}$. The diameter of the Cayley graph of this is studied by Larsen in~\cite{la} to resolve certain conjectures of Luboztky.
\end{example}

\subsection{Two-generation of finite simple groups}
In 1930, Brahana~\cite{br} noticed that several known finite simple groups are $2$-generated. Following that, Albert and Thompson (see~\cite{at}) proved the same for projective linear groups. Eventually, Steinberg in~\cite{st1}, proved that all the finite groups of Lie type are two-generated. There are several other results, case-by-case, on this subject, which we do not go into at the moment. The need of doing this case-by-case is because of the more general $(2,3)$-generation problem explained in the following section. Steinberg, following the general Chevalley-Steinberg construction, gave explicit generators in each case. For the purpose of demonstration we give an example here.
\begin{example}
Let $q$ be an odd prime power and $n\geq 2$. Fix a cyclic generator of the finite field, say, $\mathbb F_q^*=\langle \zeta \rangle$. Steinberg proved that 
$SL_n(\mathbb F_q) = \left \langle r, s \right \rangle$ where 
$$r= \begin{pmatrix} \zeta & 0 & 0\\ 0&\zeta^{-1}&0 \\ 0&0&I_{n-2}\end{pmatrix}, s= \begin{pmatrix} 1 & 0& \cdots&0& (-1)^{n-1} \\ 1& 0 &\cdots &0&0\\ &\ddots&\ddots&&\vdots \\ &&\ddots&\ddots&\vdots\\ 0&0&\cdots&1&0 \end{pmatrix}. $$
Notice that $s=x_{1,2}(1)n$, where 
$$n=\begin{pmatrix} 0& 0& \cdots&0& (-1)^{n-1} \\ 1& 0 &\cdots &0&0\\ &\ddots&\ddots&&\vdots \\ &&\ddots&\ddots&\vdots\\ 0&0&\cdots&1&0 \end{pmatrix}.$$
Let $H=\langle r, s \rangle$ be the subgroup of $SL_n(\mathbb F_q)$ generated by $ r$ and $s$. We need to prove that $H=SL_n(\mathbb F_q)$. The steps are as follows:
\begin{enumerate}
\item First we compute the commutator $ [s rs^{-1}, r] = \begin{pmatrix} 1 & -(\zeta-\zeta^{-1})^2 &0 \\ 0&1&0 \\ 0&0&I_{n-2}\end{pmatrix}$, which is in $H$. Clearly $\zeta-\zeta^{-1}\neq 0$, thus we get a non-trivial unipotent $x_{1,2}(t)$ in $H$ where $t= -(\zeta-\zeta^{-1})^2 $. By taking powers of $r$, and doing similar computations, we prove that $x_{1,2}(t)\in H$ for all $t\in \mathbb F_q$.
\item Thus, $x_{1,2}(-1).s = n$ is in $H$.
\item Now compute $nx_{1,2}(t)n^{-1}=x_{2,3}(t)$, and inductively check that $nx_{i-1,i}(t)n^{-1}=x_{i,i+1}(t)$ for all $1\leq i \leq n-1$. 
\item Now, $n x_{n-1,n}(t) n^{-1}= x_{2,1}(t)$ is in $H$, and inductively check that $nx_{i,i-1}(t)n^{-1}=x_{i+1,i}(t)$  for all $2\leq i \leq n-1$. 
\end{enumerate}
The proof is complete, combined with the fact in Example~\ref{sl-gen}. 
\end{example}

\begin{example}
Suppose $q$ is an odd prime power. The group $Sp_6(\mathbb F_q)$ is generated by the two of its elements $r$ and $s$. The element $r= diag(\zeta, \zeta^{-1}, 1, \zeta^{-1}, \zeta, 1)$, where $\zeta$ is a cyclic generators of $\mathbb F_q^*$, and 
$$s=x_{1,2}(1).w = \begin{pmatrix} -1 &0 &0&0&0&1\\ -1 &0&0&0&0&0\\ 0 & -1 &0&0&0&0\\ 0&0&-1&0&0&0 \\  0&0&1&-1&0&0\\ 0&0&0&0& -1& 0 \end{pmatrix},$$ where $w=\begin{pmatrix} 0&1\\ -I_{5}&0 \end{pmatrix}$. Let $H=\langle r, s\rangle$ be the subgroup of $Sp_6(\mathbb F_q)$, generated by $r$ and $s$. To prove equality, we prove that all the Chevalley generators are in the subgroup $H$, and hence $H= Sp_6(\mathbb F_q)$. The steps to prove this are as follows:
\begin{enumerate}
\item $x_{1,2}(\delta)$, where $\delta= (\zeta^2-1)(\zeta^{-1}-1)$, is the commutator $[s r s^{-1}, r]$, which is in $H$. By varying the powers of $r$, and multiplying the $x_{1,2}(d)$ thus obtained, we prove that $x_{1,2}(t)$ in $H$, for all $t\in \mathbb F_q$.
\item Now $w=x_{1,2}(-1).s$ is in the group $H$. 
\item $wx_{1,2}(t)w^{-1} = x_{2,3}(t)$ is in $H$. Further, $wx_{2,3}(t)w^{-1} = x_{3,-1}(t)$ is in $H$.
\item $wx_{3,-1}(t)w^{-1}=x_{2,1}(-t)$, $wx_{2,1}(t)w^{-1}=x_{3,2}(t)$, and $wx_{3,2}(t)w^{-1}=x_{-1,3}(t)$.
\end{enumerate}
Combined with the fact in Example~\ref{sp-gen}, the proof is complete,  
\end{example}

At this point, we refer the reader to survey articles by Di Martino and Tamburini~\cite{dt} and Wilson~\cite{wij} on this subject. The work on $(2,3)$-generation, discussed in the next section, naturally contains information about the two-generation problem. In our discussion we have not taken in account how many relations are required. We fast forward to the recent breakthrough in 2011 on this problem, by Guralnick, Kantor, Kassabov and Lubotzky (see Theorem A~\cite{gkkl}).
\begin{theorem}
Every finite quasi-simple (the groups which are perfect and simple modulo its center) group of Lie type, except the Ree groups ${}^2G_2(3^{2e+1})$, have presentations with $2$ generators and $51$ relations.
\end{theorem}
\noindent We urge the readers to take a moment to soak into this glorious theorem by reminding themselves that the size of finite groups of Lie type, such as the classical groups FSG3a, vary with two variables $n$ and $q$. 

\subsection{Two generation over $\mathbb Z$}
In~\cite{gt}, Gow and Tamburini proved that the group $SL_n(\mathbb Z)$, for $n\neq 4$, is generated by two matrices $x=I+\sum_{i=1}^n e_{i,i+1}$ and $y=I+\sum_{i=1}^n e_{i,i-1}$. Kassabov in~\cite{ka} has proved that the matrix ring $M_n(\mathbb Z)$, for $n\geq 2$, has a ring presentation by $2$ generators and $3$ relations as follows:
$$M_n(\mathbb Z) = \langle x, y \mid x^n, y^n, xy+y^{n-1}x^{n-1}-1 \rangle, $$
given by associating $x=\sum_{i=1}^n e_{i,i+1}$ and $y=\sum_{i=1}^n e_{i,i-1}$.
There is a lot of literature on generating dense subgroups in an arithmetic group, but we won't delve in that direction. We refer the reader to a survey article by Tamburini~\cite{ta} on this subject, and a paper by Vsemirnov~\cite{vs} for further reading.

\subsection{Standard generators used in MAGMA}
The generators for finite groups of Lie type used in MAGMA are not Chevalley generators or Steinberg $2$-generators mentioned above. These are, what is called ``Standard generators'' used for the classical groups. These were originally proposed by Costi~\cite{cos}, and later adopted by Leedham-Green and O'Brien~\cite{lo1,lo2}. These are at most $8$ in number for all classical groups.  

\section{The $(2,3)$-generation problem}
Let $G$ be a group.
\begin{definition}
A two generated group $G$ is is said to be $(2,3)$-generated if one of the generators is of order $2$ and other of order $3$. That is $G=\langle r,s \rangle$ where $r^2=1$ and $s^3=1$.
\end{definition} 
\noindent The modular group $PSL_2(\mathbb Z)$ (discussed in the Example~\ref{modular-group}) is a free product of $C_2$ and $C_3$ thus freest possible $(2,3)$-group. Thus, any $(2,3)$-generated group is a quotient of the modular group $PSL_2(\mathbb Z)$. Since finite simple groups have an element of order $2$, the question here is to classify which finite simple groups are $(2,3)$-generated, i.e., which finite simple groups are quotients of $PSL_2(\mathbb Z)$. Historically, the known answers indicated that almost all of them (see below for further details) are quotients, and thus it led to the belief that $PSL_2(\mathbb Z)$ is the mother of almost all finite groups.

Several results have been obtained for the $(2,3)$-generation of groups of Lie type: by Tamburini and Wilson~\cite{tw1, tw2} for classical groups when $n$ is large enough, Di Martino  and Vavilov~\cite{dv1,dv2} for $SL(n,q)$, Lubeck and Malle~\cite{lm} for exceptional groups, Pellegrini~\cite{pe}, Malle, Saxl and Weigel~\cite{msw} and others. In most of these results an explicit generating set is exhibited and often they prove that these groups are $(2,3)$-generated. However, there are some groups of Lie type which are not $(2,3)$-generated (see~\cite{vs}). We refer the reader to the survey article on this subject by Pellegrini and Tamburini~\cite{pt} and by Burness~\cite{bu}.

A more general problem is the $(a,b)$-generation problem, where $a$ and $b$ are given positive integers. The problem asks what are all possible finite simple groups which are quotients of the group $C_a\star C_b$. King~\cite{ki1} has proved that every non-abelian finite simple group is $(2,p)$-generated, where $p$ is a prime. We remark that in this result $p$ depends on the group. More generally, King is~\cite{ki2} has classified which finite simple groups are $(a,b)$-generated.

Clearly, for a finite simple group G to be $(2,3)$-generated, it is necessary that $G$ should possess an element of order $2$, and an element of order $3$. It is well known that non-abelian simple groups contain elements of order $2$. However, Suzuki groups, (see FSG4, Suz) do not contain element of order $3$ and hence are not $(2,3)$-generated. These are the only non-abelian simple groups which do not contain an element of order $3$. We can ask further question that if a finite simple group contains an element of order $3$, is it always $(2,3)$-generated? The answer is no! Among the classical groups $PSp(4,2^k)$ and $PSp(4,3^k)$ are not $(2,3)$-generated (see~\cite{ls}, Theorem 1.6), although they contain elements of order $3$.  For the recent update on this problem please see the survey article by~\cite{bu}, for example, Theorem 2.4.

\subsection{Probabilistic two generation}
The subject of probabilistic group theory is a less explored domain. We refer a reader to the beautiful survey article by Dixon~\cite{di} to get a feel of the subject. To give some idea of the questions dealt with in this subject, we mention some recent results proved by Liebeck and Shalev~\cite{ls}. Let $G$ be a finite simple group. The question is to understand if the probability that two randomly chosen elements of $G$ generate $G$, tends to $1$ as $|G|\rightarrow \infty$. This question is further refined with stricter conditions that one of the two randomly chosen elements is an involution. This is Conjecture 2 in~\cite{ls}, which they prove for classical and alternating groups.  The Conjecture 3 in~\cite{ls}, due to Di Martino, Vavilov, Wilson and others, is that all finite simple classical groups (with some small exception in low rank and low characteristic) are $(2,3)$-generated. Liebeck and Shalev (see Theorem 1.4 and 15 in~\cite{ls}) proved that if $G$ is a finite simple classical or alternating group, except $PSp_4(q)$, then the probability that a randomly chosen order $2$ element and a randomly chosen order $3$ element of $G$ generates $G$, tends to $1$ as $|G|\rightarrow \infty$. Thus they establish that $PSp_4(q)$ is not $(2,3)$-generated, and almost all finite simple classical groups are $(2,3)$-generated.

\subsection{$(2,3,7)$-generation}
A group is said to be $(2,3,7)$-generated if it is generated by two elements $x,y$, where $x^2=1, y^3=1$, and the product $(xy)^7=1$. Thus such groups are a quotients of the $(2,3,7)$-triangular group $\Delta$ (see Example~\ref{hurwitz}). The question is to determine all Hurwitz groups. In particular, determine which finite simple groups are $(2,3,7)$-generated, i.e., which finite simple groups are Hurwitz groups. We refer the reader to survey articles by Conder~\cite{co1,co2} on this subject.



\end{document}